\documentclass[12pt,amscd]{amsart}
\usepackage[all]{xy}
\usepackage{graphicx}
\usepackage{amsmath,amsxtra,amssymb,latexsym, amscd,amsthm, enumerate}

\usepackage{tikz}
\usepackage{ytableau}

 \usepackage{indentfirst}
\usepackage[mathscr]{eucal}
  \usepackage[pagebackref=true]{hyperref}

\newtheorem{thm}{Theorem}[section]
\newtheorem{cor}[thm]{Corollary}
\newtheorem{lem}[thm]{Lemma}
\newtheorem{prop}[thm]{Proposition}
\theoremstyle{definition}
\newtheorem{defn}[thm]{Definition}
\newtheorem{exam}[thm]{Example}

\newtheorem{rem}[thm]{Remark}

\numberwithin{equation}{section}

\DeclareMathOperator{\pd}{pd}

\begin{document}

\title{Betti numbers of powers of Ferrers ideals and edge ideals of threshold graphs}

\author{Do Trong Hoang}
\address{Faculty of  Mathematics and Informatics, Hanoi University of Science and Technology, 1 Dai Co Viet, Bach Mai, Hanoi, Vietnam.}

\email{hoang.dotrong@hust.edu.vn}

\author{Thanh Vu}
\address{Institute of Mathematics, VAST, 18 Hoang Quoc Viet, Hanoi, Vietnam}
\email{vuqthanh@gmail.com}

\subjclass[2020]{13D02, 05E40, 13F55}
\keywords{Betti numbers; Ferrers ideals; Ferrers graphs; threshold graphs; specialization}

\date{}

\commby{}
  \begin{abstract}
  We  provide explicit formulas for   the Betti numbers of all powers of Ferrers ideals and    edge ideals  of   threshold graphs.  
\end{abstract}

\maketitle

\section{Introduction}
Let $\lambda = (\lambda_1, \ldots, \lambda_n)$ be a partition, meaning that $\lambda_1 \ge \lambda_2 \ge \cdots \ge \lambda_n \ge 1$, and let $m = \lambda_1$. Let $R = k[x_1, \ldots, x_n, y_1, \ldots, y_m]$ be a polynomial ring over the field $k$. The Ferrers ideal associated with the partition $\lambda = (\lambda_1, \ldots, \lambda_n)$  is given by 
$$I_\lambda = (x_iy_j \mid 1 \le i \le n, 1 \le j \le \lambda_i).$$ 
For each partition $\lambda$, Corso and Nagel \cite{CN1} introduced and studied the Ferrers ideal $I_\lambda$ and described its minimal free resolution. In this work, we provide a formula for the Betti numbers of all powers of $I_\lambda$. Note that there are very few classes of ideals for which the Betti numbers of all powers are known, including maximal ideals, path ideals of paths \cite{BCV1}, and certain path ideals of cycles \cite{BCV2}. To achieve this, we provide an ordering of the minimal generators of $I_\lambda^t$ and show that this ordering yields a linear quotient, from which the formula for the Betti numbers follows. We now introduce these concepts in more detail to state our main results.

Each finitely generated graded $R$-module, and in particular $R/I_\lambda$, admits a minimal graded free resolution
\[
    0 \leftarrow R/I_\lambda \leftarrow \bigoplus_{j} R(-j)^{\beta_{0,j}(R/I_\lambda)} \leftarrow \cdots \leftarrow \bigoplus_{j} R(-j)^{\beta_{p,j}(R/I_\lambda)} \leftarrow 0,
\]
where $R(-j)$ denotes the free $R$-module obtained by shifting the degrees of $R$ by $j$. The number $\beta_{i,j}(R/I_\lambda)$ is the \emph{$(i,j)$-th graded Betti number} of $R/I_\lambda$. Herzog and Takayama \cite{HT} introduced the notion of ideals having linear quotients and proved that, in this case, the resolution can be constructed using iterated mapping cones. In particular, if all generators have the same degree, the ideal has a linear free resolution. We now describe the set of minimal generators of $I_\lambda^t$. 

For each $t \ge 1$, let 
$$
\begin{aligned}
C(\lambda,t) = \{ &(i_1, \ldots, i_t; j_1, \ldots, j_t) \mid 1 \le i_1 \le \cdots \le i_t \le n, \\
&j_1 \ge j_2 \ge \cdots \ge j_t, j_a \le \lambda_{i_a} \text{ for all } a = 1, \ldots, t \}.
\end{aligned}
$$
Then we have:

\begin{thm}\label{thm_linearquotient}
    The minimal generators of $I_\lambda^t$ are given by $x_{i_1} \cdots x_{i_t} y_{j_1} \cdots y_{j_t}$ for $(i_1, \ldots, i_t; j_1, \ldots, j_t) \in C(\lambda,t)$. Furthermore, when these generators are ordered lexicographically with $x_1 > x_2 > \cdots > x_n > y_1 > \cdots > y_m$, $I_\lambda^t$ has a linear quotient with respect to this order.
\end{thm}

From this, we deduce the following formula for the Betti numbers:

\begin{thm}\label{thm_betti} 
    Let $\lambda = (\lambda_1, \ldots, \lambda_n)$ be a partition with $\lambda_1 = m$. Then for all $p$, we have 
    $$\beta_p^R(I_\lambda^t) = \sum_{(i_1, \ldots, i_t; j_1, \ldots, j_t) \in C(\lambda,t)} \binom{i_t + j_1 - 2}{p}.$$
\end{thm}

The linear quotient property for powers of Ferrers graphs under lexicographic order extends to the class of threshold graphs. This allows us to deduce the Betti numbers of powers of edge ideals of threshold graphs. The remainder of this paper is devoted to proving these formulas and discussing several of their applications.

\section{Linear quotients and Betti numbers}
First, we establish the description of the minimal generating set.

\begin{lem}\label{lem:gens} 
Let $\lambda = (\lambda_1,\ldots, \lambda_n)$ be a partition with $\lambda_1 = m$. Then $I_\lambda^t$ is minimally generated by $x_{i_1}\cdots x_{i_t} y_{j_1} \cdots y_{j_t}$ for $(i_1, \ldots, i_t; j_1, \ldots, j_t) \in C(\lambda,t)$.
\end{lem}

\begin{proof} 
We have $x_{i_a} y_{j_a} \in I_\lambda$ for all $a = 1,\ldots,t$; hence, $u = x_{i_1} \cdots x_{i_t} y_{j_1} \cdots y_{j_t} \in I_\lambda^t$. Conversely, assume that $f = (x_{\alpha_1}y_{\beta_1}) \cdots (x_{\alpha_t}y_{\beta_t}) \in I_\lambda^t$ is a generator. We may assume without loss of generality that $\alpha_1 \le \cdots \le \alpha_t$. We can repeatedly apply the exchange relation to obtain the canonical representation. Indeed, suppose that we have two factors $x_i y_a$ and $x_k y_b$ with $i \le k$ and $a < b$. Then, by definition, $x_i y_b \in I_\lambda$ and $x_k y_a \in I_\lambda$. Thus, we can perform the exchange 
\[
(x_i y_a) (x_k y_b) = (x_i y_b) (x_k y_a).
\]
Applying this operation finitely many times rearranges the $y$-indices into non-increasing order without altering the $x$-indices, which completes the proof.
\end{proof}

\begin{defn}
Let $I=(u_1,\ldots,u_s)$ be a monomial ideal with ordered minimal generators $u_1,\ldots,u_s$. We say that $I$ has \emph{linear quotients} with respect to this order if, for every $j=2,\ldots,s$, the colon ideal
\[
(u_1,\ldots,u_{j-1}):u_j
\]
is generated by variables.
\end{defn}

Before proving the result on linear quotients, let us consider an example for clarity.

\begin{exam} 
Let $\lambda = (3,2)$. Then the minimal generators of $I_\lambda^2$ are ordered lexicographically as follows:
\[
\begin{aligned}
x_1^2 y_1^2 &> x_1^2y_1y_2 > x_1^2 y_1 y_3 > x_1^2y_2^2 > x_1^2y_2y_3 > x_1^2y_3^2 > x_1x_2y_1^2 \\
&> x_1x_2y_1y_2 > x_1x_2y_1y_3 > x_1x_2y_2^2 > x_1x_2y_2y_3 > x_2^2 y_1^2 > x_2^2 y_1y_2 > x_2^2 y_2^2.
\end{aligned}
\]
With this ordering, we can verify the linear quotients directly; for instance,
\[
(x_1^2y_1^2, x_1^2y_1y_2) : (x_1^2y_1y_3) = (y_1, y_2).
\]
Note that while the canonical tuples $(i_1, \ldots, i_t; j_1, \ldots, j_t) \in C(\lambda, t)$ write the $y$-indices in non-increasing order ($j_1 \ge \cdots \ge j_t$), writing the monomial factors with $y$-indices appearing in non-decreasing order aligns naturally with the lexicographic ordering.
\end{exam}

\begin{lem}\label{lem_linearquotients} 
Let $u = x_{i_1} \cdots x_{i_t} y_{j_1} \cdots y_{j_t}$ be a minimal generator of $I_\lambda^t$ corresponding to the tuple $(i_1, \ldots, i_t; j_1, \ldots, j_t) \in C(\lambda,t)$. Let $Q_u$ denote the ideal generated by all minimal generators preceding $u$ in the lexicographic order. Setting $r = i_1$ and $c = j_1$, we have 
\[
Q_u : u = (x_1, \ldots, x_{r-1}, y_1, \ldots, y_{c-1}).
\]
\end{lem}

\begin{proof} 
Let $r = i_1$ and $c = j_1$. 
First, we prove the inclusion $(x_1, \ldots, x_{r-1}, y_1, \ldots, y_{c-1}) \subseteq Q_u : u$. 
For any $a < r$, replacing $x_r$ with $x_a$ yields a valid minimal generator $v$ in $I_\lambda^t$ (after sorting the indices into $C(\lambda,t)$) that strictly precedes $u$ in the lexicographic order, meaning $x_a \in Q_u : u$. 
Similarly, for any $b < c$, replacing $y_c$ with $y_b$ yields a valid minimal generator $v$ preceding $u$, so $y_b \in Q_u : u$. 

Conversely, to show that $Q_u : u \subseteq (x_1, \ldots, x_{r-1}, y_1, \ldots, y_{c-1})$, recall the variable ordering $x_1 > \cdots > x_n > y_1 > \cdots > y_m$. In the lexicographic order, any preceding generator $v <_{\text{lex}} u$ must introduce a strictly larger variable at the earliest possible position where $v$ and $u$ differ. Since $u$ has leading variable components determined by $i_1 = r$ and $j_1 = c$, any variable dividing $v / \gcd(u,v)$ that places $v$ before $u$ must belong to the set $\{x_1, \ldots, x_{r-1}\}$ or $\{y_1, \ldots, y_{c-1}\}$. Consequently, every generator of $Q_u : u$ is divided by at least one variable from this set, establishing the reverse inclusion and completing the proof.
\end{proof}

\begin{proof}[Proof of Theorem \ref{thm_linearquotient}] 
The conclusion follows directly from Lemma \ref{lem:gens} and Lemma \ref{lem_linearquotients}.    
\end{proof}

\begin{proof}[Proof of Theorem \ref{thm_betti}] 
The assertion follows by induction from Lemma \ref{lem_linearquotients} together with \cite[Theorem 2.6]{SV}.    
\end{proof}

\begin{exam}
Let $\lambda=(3,2)$ and $t=4$, so that $R=k[x_1,x_2,y_1,y_2,y_3]$. Each tuple in $C(\lambda,4)$ can be written uniquely as
\[
\bigl((\underbrace{1,\ldots,1}_{a},\underbrace{2,\ldots,2}_{4-a});
      (\underbrace{3,\ldots,3}_{q},\underbrace{2,\ldots,2}_{r},
       \underbrace{1,\ldots,1}_{4-q-r})\bigr),
\qquad 0\le q\le a\le 4,\quad 0\le r\le 4-q.
\]
Here, $a$ counts the occurrences of $1$ among the $i$-indices, while $q$ and $r$ count the occurrences of $3$ and $2$ among the $j$-indices. Therefore,
\[
|C(\lambda,4)|=\sum_{a=0}^4\sum_{q=0}^a(5-q)=55.
\]
The choices of $(a,q,r)$ can be grouped according to the value of $i_4+j_1-2$ as follows:
\[
\begin{array}{c|c|c}
\text{Parameter choices} & \text{Number of tuples} & i_4+j_1-2 \\
\hline
0\le a\le3,\ q=0,\ r=0 & 4 & 1 \\
0\le a\le3,\ q=0,\ 1\le r\le4 & 16 & 2 \\
1\le q\le a\le3,\ 0\le r\le4-q & 20 & 3 \\
a=4,\ q=0,\ r=0 & 1 & 0 \\
a=4,\ q=0,\ 1\le r\le4 & 4 & 1 \\
a=4,\ 1\le q\le4,\ 0\le r\le4-q & 10 & 2
\end{array}
\]
Thus, the exponent values $0,1,2,3$ occur respectively $1,8,26,20$ times. By Theorem \ref{thm_betti},
\[
\beta_p^R(I_{\lambda}^4)
=\binom{0}{p}+8\binom{1}{p}+26\binom{2}{p}+20\binom{3}{p}.
\]
Evaluating this expression gives $\beta_0=55$, $\beta_1=120$, $\beta_2=86$, $\beta_3=20$, and $\beta_p=0$ for $p\ge4$. Hence, the Betti sequence is $(55,120,86,20)$.
\end{exam}

\begin{rem} 
It is well known that $I_\lambda^t$ has linear resolutions for all $t \ge 1$ (see, e.g., \cite[Theorem 3.2]{HHZ}). The depth of powers of $R/I_\lambda^t$ is also well known (see, e.g., \cite[Proposition 2.18]{HV}). Basser, Diethorn, Miranda, and Stinson-Maas \cite[Theorem 3.11]{BDMS} show that powers of edge ideals of cochordal graphs have linear quotients. Theorem \ref{thm_linearquotient} provides a specific ordering that yields linear quotients for $I_\lambda^t$, resulting in a more explicit formula for the Betti numbers.
\end{rem}

\section{Betti numbers of powers of threshold graphs}

Threshold graphs are a special class of cochordal graphs, where Ferrers graphs correspond precisely to bipartite threshold graphs. They have numerous applications (see, e.g., \cite{MP}). Corso and Nagel \cite{CN2} showed that each threshold graph can be obtained by specializing a Ferrers graph. More precisely, by \cite[Corollary 4.6]{CN2}, every threshold graph can be characterized by a partition $\lambda = (\lambda_1, \ldots, \lambda_n)$ with $\lambda_1 \ge \lambda_2 \ge \cdots \ge \lambda_n > n$ where $\{i,j\}$ with $i < j$ is an edge of the threshold graph $T_\lambda$ if and only if $j \le \lambda_i$. Let $S = k[z_1, \ldots, z_m]$ where $m = \lambda_1$. The corresponding edge ideal of $T_\lambda$ is given by 
\[
I(T_\lambda) = (z_i z_j \mid i < j \le \lambda_i).
\]
For a monomial $u = z_1^{a_1} \cdots z_t^{a_t}$ with $a_t > 0$, we denote by $\operatorname{last}(u) = t$ the largest index of a variable dividing $u$.

\begin{thm}\label{thm_linear_quotient_threshold} 
The minimal generators of $I(T_\lambda)^t$, ordered lexicographically with $z_1 > z_2 > \cdots > z_m$, form a linear quotient order for $I(T_\lambda)^t$. Furthermore, for a minimal generator $u = z_1^{a_1} \cdots z_m^{a_m}$ of $I(T_\lambda)^t$, we have 
\[
Q_u : u = (z_r \mid a_r < t \text{ and } r < \operatorname{last}(u)),
\]
where $Q_u$ is the ideal generated by all minimal generators preceding $u$ in this order.
\end{thm}

\begin{proof} 
Let $q = \operatorname{last}(u)$ and let $L = (z_r \mid a_r < t \text{ and } r < q)$. 
First, we claim that $L \subseteq Q_u : u$. Assume that $r < q$ and $a_r < t$. Choose a factorization of $u$ into a product of $t$ edges, $u = e_1 \cdots e_t$. Since $q \mid u$, there exists an edge $e_j$ such that $z_q$ divides $e_j$. Hence, we may write $e_j = z_p z_q$ with $p < q$. Since $e_j \in I(T_\lambda)$, we have $q \le \lambda_p$. There are two cases:

\textbf{Case 1: $p \neq r$.} We claim that $z_p z_r \in I(T_\lambda)$. Indeed, if $p < r$, then $p < r < q \le \lambda_p$, hence $z_p z_r \in I(T_\lambda)$. If $r < p$, then $p < q \le \lambda_p \le \lambda_r$, hence $z_p z_r \in I(T_\lambda)$. In either case, we have $w = e_1 \cdots e_{j-1} (z_p z_r) e_{j+1} \cdots e_t \in I(T_\lambda)^t$. Furthermore, $w / \gcd(w,u) = z_r$. Thus, $z_r \in Q_u : u$.

\textbf{Case 2: $p = r$.} Since $a_r < t$, there exists an edge $e_i$ such that $z_r$ does not divide $e_i$. Write $e_i = z_a z_b$. Since $q = \operatorname{last}(u)$, we have $a, b \le q$. We claim that $z_r z_a \in I(T_\lambda)$. Indeed, since $z_r z_q \in I(T_\lambda)$, we have $r < q \le \lambda_r$. Now, if $r < a$, then $r < a < q \le \lambda_r$, hence $z_r z_a \in I(T_\lambda)$. If $a < r$, then $a < r < q \le \lambda_r \le \lambda_a$, hence $z_r z_a \in I(T_\lambda)$. Similarly, $z_r z_b \in I(T_\lambda)$. Replacing $e_i e_j$ by $(z_r z_a)(z_r z_b)$ yields a new element $w \in I(T_\lambda)^t$ such that $w / \gcd(w,u) = z_r$, which implies $z_r \in Q_u : u$.

Now, let $v$ be any generator preceding $u$ in the lexicographic order. Write $v = z_1^{b_1} \cdots z_m^{b_m}$. Let $r$ be the smallest index such that $b_r \neq a_r$. By definition of the lexicographic order, we have $b_r > a_r$ and $b_i = a_i$ for all $i < r$. Since $v$ is a minimal generator of $I(T_\lambda)^t$, we must have $b_r \le t$, which implies $a_r < b_r \le t$. Furthermore, $\deg(v) = \deg(u)$, so there must exist an index $j > r$ such that $b_j < a_j$. In particular, $r < j \le \operatorname{last}(u) = q$. Consequently, $v / \gcd(u,v)$ is divisible by $z_r$, showing that every generator of $Q_u : u$ is divisible by a variable $z_r$ with $r < q$ and $a_r < t$. The conclusion follows.
\end{proof}

From Theorem \ref{thm_linear_quotient_threshold}, we see that to find the Betti numbers of the powers $I(T_\lambda)^t$, we need to determine formulas for two key invariants. 
First, let 
\[
\operatorname{Gens}_t(\lambda, q) = \{ u \mid u \text{ is a minimal generator of } I(T_\lambda)^t \text{ with } \operatorname{last}(u) = q \}
\]
and 
\[
\operatorname{Gens}^*_t(\lambda, q) = \{ u = z_1^{a_1} \cdots z_q^{a_q} \in \operatorname{Gens}_t(\lambda, q) \mid a_i = t \text{ for some } i < q \}.
\]
We then define the cardinalities:
\begin{align*}
    a_t(\lambda; q) &= |\operatorname{Gens}_t(\lambda, q)|, \\
    b_t(\lambda; q) &= |\operatorname{Gens}^*_t(\lambda, q)|.
\end{align*}

\begin{thm}\label{thm_Betti_threshold} 
Let $T_\lambda$ be the threshold graph corresponding to a partition $\lambda$ with $\lambda_1 \ge \cdots \ge \lambda_n \ge n$. Then 
\[
\beta_p(I(T_\lambda)^t) = \sum_{q=2}^m \left[ a_t(\lambda; q) \binom{q-1}{p} - b_t(\lambda; q) \binom{q-2}{p-1} \right].
\]
\end{thm}

\begin{proof}
By Theorem \ref{thm_linear_quotient_threshold} and \cite[Theorem 2.6]{SV}, we have 
\[
\beta_p(I(T_\lambda)^t) = \sum_{u \in G(I_\lambda^t)} \binom{|C(u)|}{p},
\]
where $C(u)$ is the set of variables in the quotient ideal $Q_u : u$. For any minimal generator $u$ with $\operatorname{last}(u) = q$, the size of $C(u)$ satisfies
\[
|C(u)| = 
\begin{cases} 
q - 2 & \text{if } u \in \operatorname{Gens}_t^*(\lambda, q), \\ 
q - 1 & \text{otherwise}. 
\end{cases}
\]
The conclusion follows.
\end{proof}

\begin{cor}\label{cor_pd_threshold} 
Let $T_\lambda$ be the threshold graph corresponding to a partition $\lambda$ with $\lambda_1 = m \ge \cdots \ge \lambda_n \ge n$. Then for all $t \ge 2$,
\[
\pd(I(T_\lambda)^t) = 
\begin{cases} 
m - 1 & \text{if } n \ge 2, \\ 
m - 2 & \text{if } n = 1. 
\end{cases}
\]
\end{cor}

\begin{proof}
If $n = 1$, then $I(T_\lambda) = z_1(z_2,\ldots,z_m)$, and its powers have projective dimension $m - 2$. We may therefore assume that $n \ge 2$ and $t \ge 2$, which implies $\lambda_1 = m > n \ge 2$. By the Betti number formula in Theorem \ref{thm_Betti_threshold}, it suffices to show that there exists a minimal generator $u$ of $I(T_\lambda)^t$ such that $|C(u)| = m - 1$. To this end, consider the monomial $u = (z_1 z_m)(z_2 z_3)^{t-1}$. Then $u \in \text{Gens}_t(\lambda, m)$ but $u \notin \text{Gens}_t^*(\lambda, m)$. Consequently, $|C(u)| = m - 1$. The conclusion follows.
\end{proof}

\begin{lem}\label{lem_b_t} 
Let $\rho_\lambda(q) = |\{ 1 \le r < q \mid \lambda_r \ge q \}$. Then 
\[
b_t(\lambda; q) = \rho_\lambda(q) \binom{q + t - 3}{t - 1}.
\]
\end{lem}

\begin{proof} 
Let $u = z_1^{a_1} \cdots z_q^{a_q}$ be a minimal generator of $I(T_\lambda)^t$ such that $a_r = t$ for some index $r < q$. This condition implies that the edge $z_r z_q$ belongs to $I(T_\lambda)$, which is equivalent to $\lambda_r \ge q$. Under this condition, we have $z_r z_i \in I(T_\lambda)$ for all $i \in \{1, \ldots, q\} \setminus \{r\}$.

There are $\rho_\lambda(q)$ possible choices for the index $r < q$ such that $\lambda_r \ge q$. For a fixed $r$, the condition $a_r = t$ means that $z_r^t$ accounts for $t$ factors of $u$. The remaining part is a monomial of degree $t$ in the variables $\{z_1, \ldots, z_q\} \setminus \{z_r\}$. The number of such monomials is given by $\binom{q + t - 3}{t - 1}$. Summing over all $\rho_\lambda(q)$ possible choices for $r$, we obtain
\[
b_t(\lambda; q) = \rho_\lambda(q) \binom{q + t - 3}{t - 1}.
\]
This completes the proof.
\end{proof}

To compute $a_t(\lambda, q)$, we introduce the following notation. For a partition $\lambda$ and an integer $q$, we define the truncation partition $\lambda^{[q]}$ to be the partition given by $\lambda^{[q]}_i = \min \{q, \lambda_i\}$, and we let $\bar{\lambda}$ denote the partition $(\lambda_2 - 1, \ldots, \lambda_n - 1)$. For any $r, s \ge 0$, we set 
\[
A(\lambda; r, s) = \dim_k (I(T_\lambda)^t)_{2r + s}.
\]

In particular, $A(\lambda; t, 0) = \mu(I (T_\lambda)^t)$, which is the number of minimal generators of $I (T_\lambda)^t$.

First, we have the following lemma:
\begin{lem}\label{lem_a_1} 
With the notation above, we have 
\[
a_t(\lambda, q) = A(\lambda^{[q]}; t, 0) - A(\lambda^{[q-1]}; t, 0).
\]
\end{lem}
\begin{proof}
Any generator of $I_{\lambda^{[q]}}^t$ that does not belong to $I_{\lambda^{[q-1]}}^t$ is divisible by $z_q$. The conclusion then follows directly from the definitions.
\end{proof}

\begin{lem}\label{lem_a_2} 
Assume that $2 \le q \le n+1$. Then 
\[
a_t(\lambda; q) = \binom{2t + q - 2}{q - 1} - (q - 1)\binom{t + q - 3}{q - 1} - \binom{t + q - 2}{q - 1}.
\]
\end{lem}

\begin{proof} 
When $2 \le q \le n + 1$, the truncation $\lambda^{[q]}$ gives rise to a complete graph on $q$ vertices. Thus, $a_t(\lambda; q)$ equals the number of tuples of nonnegative integers $(a_1, \ldots, a_q)$ satisfying $a_i \le t$ for all $i$, $a_q \ge 1$, and $\sum_{i=1}^q a_i = 2t$. 

Setting $b_q = a_q - 1$, this condition corresponds to tuples of nonnegative integers $(a_1, \ldots, a_{q-1}, b_q)$ such that $\sum_{i=1}^{q-1} a_i + b_q = 2t - 1$, with the constraints $a_i \le t$ for $1 \le i \le q-1$ and $b_q \le t - 1$. 

Without these upper-bound constraints, the number of such tuples-equal to the number of degree $2t - 1$ monomials in $q$ variables-is given by the first term, $\binom{2t + q - 2}{q - 1}$. 

If any of the variables $a_i$ ($1 \le i \le q-1$) violates the condition by satisfying $a_i \ge t + 1$, we can factor out $z_i^{t+1}$, leaving a monomial of degree $t - 2$ in $q$ variables. Summing over the $q-1$ choices for such an index yields the second term, $(q - 1)\binom{t + q - 3}{q - 1}$. 

Similarly, if the boundary variable satisfies $b_q \ge t$, the remaining terms form a monomial of degree $t - 1$ in $q$ variables, giving the third term, $\binom{t + q - 2}{q - 1}$. 

Since no two upper-bound violations can occur simultaneously, we obtain the desired formula.
\end{proof}

\begin{lem}\label{lem_a_3}
Assume that $\lambda_1 = \cdots = \lambda_a = m > \lambda_{a+1}$ with $a \ge 2$. Then 
\[
A(\lambda; t, s) = \binom{m + t + s - 1}{m - 1} - \binom{m + t - 2}{m - 1} + \sum_{r=1}^t A(\bar{\lambda}; r, t + s - r).
\]
\end{lem}

\begin{proof}
For simplicity of notation, let $J = I_\lambda$ and $K = I_{\bar{\lambda}}$. Then we have $J = z_1(z_2, \ldots, z_m) + K$. Consequently, the $t$-th power of $J$ admits the expansion
\[
J^t = \sum_{r=0}^t z_1^{t-r} (z_2, \ldots, z_m)^{t-r} K^r.
\]
For any monomial $u \in (J^t)_{2t+s}$, we can write $u = z_1^d v$ for some integer $d \ge 0$. When $d \le t-1$, setting $r = t-d$, the remaining part $v$ belongs to $(z_2, \ldots, z_m)^{t-r} K^r$, which corresponds to the graded component of $K^r$ in degree $2t + s - d = t + s + r = 2r + (t + s - r)$. This yields the contribution $A(\bar{\lambda}; r, t + s - r)$. 

When $d \ge t$, the monomials arise entirely from the term $(z_2, \ldots, z_m)^t$. The number of such monomials equals the number of monomials of degree $t$ through $t+s$ in a polynomial ring in $m-1$ variables. The total number of such monomials is given by
\[
\sum_{d = t}^{t+s} \binom{d+m-2}{m-2} = \binom{m + t + s - 1}{m - 1} - \binom{m + t - 2}{m - 1}.
\]
Combining these two cases yields the desired formula. The conclusion follows.
\end{proof}

\begin{lem}\label{lem_a_4} 
Assume that $\lambda_1 = m > \lambda_2$, and let $e = m - \lambda_2$. Then 
\begin{align*}
A(\lambda; t, s) &= \binom{m + t + s - 1}{m - 1} - \binom{m + t - 2}{m - 1} \\
&\quad + \sum_{r = 1}^t \sum_{v = 0}^{t+s-r} \binom{e + t + s - r - v - 1}{e - 1} A(\bar{\lambda}; r, v).
\end{align*}
\end{lem}

\begin{proof}
As in the proof of Lemma \ref{lem_a_3}, the case where $d \ge t$ yields the standard binomial difference term. When $d < t$, we must account for the graded components of $K^r$ within the polynomial ring in variables $z_2, \ldots, z_m$. Since $K$ is supported only on the variables $z_2, \ldots, z_{\lambda_2}$, the remaining $e = m - \lambda_2$ variables lie outside the support of $K$. The inner sum over $v$ correctly accounts for the contribution of these $e$ additional variables.
\end{proof}

\begin{cor}\label{cor_Betti_complete} 
Let $K_m$ denote the complete graph on $m$ vertices. Then for all $0 \le p \le m - 1$, we have
\[
\beta_p(I(K_m)^t) = \binom{2t + m - 1}{2t + p}\binom{2t + p - 1}{p} - m\binom{t + m - 2}{m - 1}\binom{m - 1}{p}.
\]
\end{cor}

\begin{proof}
The complete graph $K_m$ corresponds to the threshold graph with $n = m - 1$ and $\lambda_1 = \cdots = \lambda_n = m$. By Theorem \ref{thm_Betti_threshold}, Lemma \ref{lem_b_t}, and Lemma \ref{lem_a_2}, we have
\[
\beta_p(I(T_\lambda)^t) = \sum_{q=2}^m \left[ a_t(\lambda; q)\binom{q-1}{p} - b_t(\lambda; q)\binom{q-2}{p-1} \right],
\]
where $a_t(\lambda; q) = \binom{2t + q - 2}{q - 1} - (q - 1)\binom{t + q - 3}{q - 1} - \binom{t + q - 2}{q - 1}$ and $b_t(\lambda; q) = (q - 1)\binom{t + q - 3}{t - 1}$. 

Setting $k = q - 1$, we can rewrite the Betti number summation as
\[
\beta_p(I(T_\lambda)^t) = \sum_{k=1}^{m-1} \binom{2t + k - 1}{k}\binom{k}{p} - \sum_{k=1}^{m-1} E_p(k),
\]
where
\[
E_p(k) = \left[ k\binom{t + k - 2}{k} + \binom{t + k - 1}{k} \right]\binom{k}{p} + k\binom{t + k - 2}{k - 1}\binom{k - 1}{p - 1}.
\]
Using the binomial identity 
\[
\binom{2t + k - 1}{k}\binom{k}{p} = \binom{2t + p - 1}{p}\binom{2t + k - 1}{k - p},
\]
we obtain
\[
\sum_{k=1}^{m-1} \binom{2t + k - 1}{k}\binom{k}{p} = \binom{2t + p - 1}{p}\binom{2t + m - 1}{2t + p} - \delta_{p,0},
\]
where $\delta_{p,0} = 1$ if $p = 0$ and $0$ otherwise. Furthermore, using the binomial identities 
\[
\binom{k}{p} = \binom{k-1}{p} + \binom{k-1}{p-1} \quad \text{and} \quad \binom{t+k-1}{k} = \binom{t+k-2}{k} + \binom{t+k-2}{k-1},
\]
we can express $E_p(k)$ as a telescoping difference $E_p(k) = F_p(k+1) - F_p(k)$, where
\[
F_p(k) = (k + 1)\binom{t + k - 1}{k}\binom{k}{p}.
\]
Summing over $k$ from $1$ to $m-1$, the terms telescope, allowing us to deduce the final formula. The conclusion follows.
\end{proof}

\begin{rem}\label{rem_arbitrary_threshold} 
Theorem \ref{thm_Betti_threshold} and Lemmas \ref{lem_b_t}--\ref{lem_a_4} allow us to compute the Betti numbers of arbitrary threshold graphs. Dung and Vu \cite{DV} associated a type sequence to any cochordal graph and gave a formula for the Betti numbers of its edge ideal in terms of this type sequence. However, we note that the type sequence alone is not sufficient to determine the Betti numbers of powers of its edge ideal. For example, $K_{2,2}$ and a triangle with a pendant edge both have the same type sequence $(2,2)$, but their second powers do not share the same Betti numbers.
\end{rem}

We do not know whether the linear quotient property with respect to the lexicographic ordering holds for powers of edge ideals of arbitrary specializations. Nonetheless, we can show that the projective dimension of powers of any arbitrary specialization is equal to $m$ if it contains a triangle. Note that if the specialization is bipartite, then its edge ideal is isomorphic to a Ferrers ideal.

\begin{prop} \label{prop_pd}
Let $\lambda = (\lambda_1, \ldots, \lambda_n)$ be a partition and $\mu = (\mu_1, \ldots, \mu_n)$ satisfy $\mu_1 \le \cdots \le \mu_n < \lambda_n$. Let $I = (z_i z_j \mid 1 \le i \le n, \mu_i < j \le \lambda_i)$ be the edge ideal of the specialization graph $G_{\lambda, \mu}$. Assume that $I$ contains a triangle, or equivalently, $\mu_1 < n$. Then $\pd(S/I^t) = m$ for all $t \ge 2$, where $m = \lambda_1$ denotes the number of vertices.
\end{prop}

\begin{proof}
First, we show that $I$ is non-bipartite if and only if $\mu_1 < n$. Indeed, if $\mu_1 \ge n$, then the underlying graph $G_{\lambda, \mu}$ is bipartite with bipartition $\{1, \ldots, n\} \cup \{\mu_1 + 1, \ldots, \lambda_n\}$. Conversely, if $\mu_1 < n$, setting $j = \mu_1 + 1$ and $k = \lambda_n$, we see that $\{1, j, k\}$ forms a triangle.

Next, we claim that the triangle $\{1, j, k\}$ dominates $G_{\lambda, \mu}$. To see this, let $v$ be any vertex:
\begin{itemize}
    \item If $v > \mu_1$, then $\{1, v\}$ is an edge by definition.
    \item If $v \le \mu_1$, then $\mu_v \le \mu_n < \lambda_n = k \le \lambda_v$, which implies that $\{v, k\}$ is an edge.
\end{itemize}

Since $G_{\lambda, \mu}$ contains a dominating triangle, it follows that $\pd(S/I^t) = m$ for all $t \ge 2$. The conclusion follows.
\end{proof}

\begin{rem}\label{rem_counterexample} 
The conclusion of Proposition \ref{prop_pd} does not hold for arbitrary cochordal graphs containing a triangle. Indeed, let 
$$I = (x_1x_2, x_1x_3, x_1x_4, x_2x_3, x_2x_4, x_3x_4, x_1y_1, x_2y_2, x_3y_3, x_4y_4)$$
be the edge ideal of the whisker graph on the complete graph $K_4$. Then $\pd(S/I^2) = 7 < 8$.
\end{rem}

\noindent\textbf{Acknowledgment.} Thanh Vu is partially supported by the Vietnam National Foundation for Science and Technology Development under the grant number NAFOSTED/101/10104/2026/1279.

\end{document}